\documentclass[12pt]{amsart}
\usepackage{calc,amssymb,amsmath,ulem}
\usepackage{alltt}
\usepackage[top=1.05in,right=1.25in,left=1.25in,bottom=1.05in]{geometry}
\RequirePackage[dvipsnames,usenames]{color}

\input{kmacros-s.sty}
\usepackage{hyperref}
\usepackage[all,cmtip]{xy}
\usepackage{verbatim}
\usepackage{mathrsfs}

\usepackage{ifthen}
\newboolean{bool.doWithPairs}
\setboolean{bool.doWithPairs}{false}

\renewcommand{\O}{\mathscr O}

\renewcommand{\to}[1][]{\xrightarrow{\ #1\ }}
\newcommand{\DBDual}[1]{{{\uline \omega}{}^{\mydot}_{\,#1}}}
\newcommand{\myH}{h}

\input{smacros2.sty}
\renewcommand{\qis}{\simeq}
\renewcommand{\O}{\mathscr O}

\begin{document}
\input{smacros1.sty}
\title{Du~Bois singularities deform}
\author{S\'andor J Kov\'acs and Karl Schwede}
\begin{abstract}
  Let $X$ be a variety and $H$ a Cartier divisor on $X$.  We prove that if $H$ has Du
  Bois (or DB) singularities, then $X$ has Du Bois singularities near $H$.  As a
  consequence, if $X \to S$ is a proper flat family over a smooth curve $S$ whose special fiber
  has Du Bois singularities, then the nearby fibers also have Du Bois singularities.
  We prove this by obtaining an injectivity theorem for certain maps of canonical
  modules.  As a consequence, we also obtain a restriction theorem for certain non-lc
  ideals.
\end{abstract}
\subjclass[2000]{14B07, 14B05, 14F17, 14F18}
\keywords{Du Bois singularities, deformation, log canonical singularities}
\address{Department of Mathematics\\ University of Washington\\ Seattle, WA, 98195,
  USA}
\email{skovacs@uw.edu}
\address{Department of Mathematics\\ The Pennsylvania State University\\ University
  Park, PA, 16802, USA}
\email{schwede@math.psu.edu}

\thanks{%
  The first named author was partially supported by in part by NSF Grant DMS
  \#0856185, and the Craig McKibben and Sarah Merner Endowed Professorship in
  Mathematics at the University of Washington}%
\thanks{%
  The second named author was partially supported by the NSF grant DMS \#1064485 and
  an NSF Postdoctoral Fellowship}

\subjclass[2010]{14B07, 14F18, 14B05} \keywords{Du Bois singularities, DB
  singularities, deformation, log canonical singularities, non-lc ideal}
\maketitle

\section{Introduction}

Du Bois singularities, or henceforth simply DB singularities, were introduced by Steenbrink in
\cite{SteenbrinkCohomologicallyInsignificant}. They may be
considered a generalization of the notion of rational singularities. The definition
and its simple consequences makes DB singularities the natural class to consider in
many situations including vanishing theorems and moduli theory. More precisely it
is important and useful that the singularities considered in these situation are
Du~Bois.  For instance, Steenbrink showed that families over smooth curves whose
fibers have DB singularities possess particularly nice properties; this maxim and
its consequences have been further explored in \cite[Section
7]{KollarKovacsLCImpliesDB}.
These applications imply that the question of whether DB singularities are invariant
under small deformations, that is whether the property of having DB singularities is
open in flat families, is very important. In this paper we settle this question in the
affirmative.

As both rational singularities \cite{KovacsDuBoisLC1} and log canonical singularities
\cite{KollarKovacsLCImpliesDB} are DB, it is interesting to note that rational
singularities are invariant under small deformations by
\cite{ElkikDeformationsOfRational}, while log canonical singularities are not unless
the total space has a $\bQ$-Cartier canonical divisor compatible with the canonical
divisors of the family members. In this latter case the statement follows from
inversion of adjunction \cite{KawakitaInversion}.

Our main result is the following:
\begin{mainthm*} [\textnormal{[Theorem \ref{thm.MainTheorem}]}]
  Let $X$ be a scheme of finite type over $\bC$ and $H$ a reduced Cartier divisor on
  $X$. If $H$ has DB singularities, then $X$ has DB singularities near $H$.
\end{mainthm*}
The openness of the Du Bois locus in proper flat families follows immediately, see Corollary
\ref{cor.FamilyDeforms}.

In \cite{IshiiSmallDeformations}, Ishii proved this result for isolated Gorenstein
singularities, and it follows for normal Gorenstein singularities from a combination
of \cite{KovacsDuBoisLC1} and \cite{KawakitaInversion}.  The first named author
claimed a proof of the same statement in general in \cite{KovacsDuBoisLC2}.  That
proof unfortunately is incomplete and only works under an additional condition.  The
problem lies in the first paragraph of the proof, namely that one may not always
reduce to the case when the non-Du Bois locus of $X$ is contained in $H = X_s$. For
additional discussion of this issue see \cite[Section 12]{KovacsSchwedeDuBoisSurvey}.

In this paper we correct that proof by showing a more general injectivity theorem,
Theorem \ref{thm.MainInjectivity}, which should be viewed as playing the same role
for Du Bois singularities that Grauert-Riemenschneider vanishing plays for rational
singularities, see Corollary \ref{cor.GRVanishingForDuBois}.  Using this injectivity,
we can follow the strategy of \cite{KovacsDuBoisLC2} and mimic Elkik's proof
\cite{ElkikDeformationsOfRational} that rational singularities deform in families to
obtain the main result.

As another corollary of this injectivity theorem, we also prove a restriction theorem
for the so-called \emph{maximal non-lc ideals} defined in
\cite{FujinoSchwedeTakagiSupplements}, at least in the case of a Gorenstein ambient
variety, see Theorem \ref{thm.RestrictionTheoremForMNLC}.
\vskip6pt
\noindent{\it{Acknowledgements:}}

The authors would like to thank the referee, Kazuma Shimomoto and Burt Totaro for pointing out typos in previous drafts of this paper.  We would also like to thank Florin Ambro for pointing out some references to us.

\section{Preliminaries on DB singularities}

Throughout this paper, all schemes are assumed to be separated and of finite type
over $\bC$, and all morphisms are defined over $\bC$.  A \emph{variety} here means a
reduced connected scheme.

We use $D^{b}_{\coherent}(X)$ to denote the bounded derived category of
$\O_X$-modules with coherent cohomology.  Given an object $C^{\mydot} \in
D^{b}_{\coherent}(X)$, its $i$th cohomology is denoted by $\myH^i(C^{\mydot})$.  For
any scheme $X$ of finite type over $\bC$, we use $\omega_X^{\mydot}$ to denote the
dualizing complex of $X$ which is defined as $\epsilon^! \bC$ where $\epsilon : X \to
\bC$ is the structure map of $X$.  We will repeatedly use Grothendieck duality in the
following form: For any proper map of schemes $f : Y \to X$, and any $C^{\mydot} \in
D^{b}_{\coherent}(Y)$, there exists a functorial quasi-isomorphism:
\[
\myR f_* \myR \sHom_{Y}^{\mydot}(C^{\mydot}, \omega_Y^{\mydot}) \qis \myR
\sHom_{X}^{\mydot}(\myR f_* C^{\mydot}, \omega_X^{\mydot}).
\]
For an introduction to derived categories and Grothendieck duality in the context
used in this paper, see \cite{HartshorneResidues}.

Recall that given a variety $X$, a \emph{resolution of singularities} $\pi : \tld X
\to X$ is a proper birational\footnote{birational = there exists a bijection of
  irreducible components with an induced isomorphism of fields of fractions.} map
from a smooth variety $\tld X$.  Given a closed subscheme $Z \subseteq X$ with
associated ideal sheaf $\sI_Z$, we say that $\pi : \tld X \to X$ is a \emph{log
  resolution of $Z \subseteq X$} if $\pi$ is a resolution of singularities and if in
addition $\pi^*\sI_Z\simeq \O_{\tld X}(G)$ where $G$ is a divisor, the exceptional
set of $\pi$, $\exc(\pi) \subseteq \tld X$, is also a divisor, and the divisor
$\exc(\pi) \cup \Supp(G)$ has simple normal crossings.  Note that resolutions of
singularities, and log resolutions, exist by \cite{HironakaResolution}.

We briefly recall some common objects used in the study of DB singularities.  For a
more extensive discussion of DB singularities, please see
\cite{KovacsSchwedeDuBoisSurvey}, \cite[Section
3.I]{HaconKovacsClassificationBook}, or \cite{PetersSteenbrinkMixedHodgeStructures}.

\begin{lemma}
\label{lem.BasicPropertiesOfDuBois}
Given a variety $X$, one may associate to $X$ an object $\DuBois{X} \in
D^{b}_{\coherent}(X)$ defined as follows: let $\pi_{\mydot} : X_{\mydot} \to X$ be a
(cubic or simplicial) hyperresolution of $X$, see \cite{GNPP, CarlsonPolyhedral,
  DeligneHodgeIII}, then
\[
\DuBois{X} := \myR {\pi_{\mydot}}_* \O_{X_{\mydot}}.
\]

This object has the following properties:
\begin{itemize}
\item[(i)] $\DuBois{X}$ is functorial with respect to morphisms of varieties, i.e.,
  given a morphism of varieties $f : Y \to X$, there is an induced morphism
  $\DuBois{X} \to \myR f_* \DuBois{Y}$.
\item[(ii)]  There is a natural morphism $\O_X \to \DuBois{X}$ compatible with (i) in the obvious way.
\item[(iii)] If in addition $X$ is proper, then the composition
\[
H^i(X^{\an}, \bC) \to H^i(X, \O_X) \to \bH^i(X, \DuBois{X})
\]
is surjective.
\end{itemize}
\end{lemma}
\begin{proof}
  See \cite{DuBoisMain} and \cite{SteenbrinkCohomologicallyInsignificant} for the
  original definitions and proofs and \cite{KovacsSchwedeDuBoisSurvey} for a survey
  on DB singularities.  Property (iii) follows directly from the $E_1$-degeneration
  of the Deligne-Du Bois variant of the Hodge-to-De Rham spectral sequence.
\end{proof}

\begin{definition}
  We say that $X$ has \emph{DB singularities} if the morphism $\O_X \to \DuBois{X}$ from
  (ii) above is a quasi-isomorphism.
\end{definition}

We also recall the following fact about DB singularities.
\begin{lemma}[(\protect{\cf~\cite[Proof of Theorem 12.8]{KollarShafarevich}})]
  \label{lem.GeneralSectionOfDuBois}
  If $X$ has DB singularities and $H$ is a general member of a base-point-free linear
  system $\delta$ on $X$, then $H$ also has DB singularities.
\end{lemma}

In this paper, we will repeatedly use the Grothendieck dual of $\DuBois{X}$.  To make
that easier we introduce the notation
\[
\DBDual{X} := \myR \sHom_{X}^{\mydot}(\DuBois{\,X}, \omega_X^{\mydot}).
\]
We will also use the fact that there exists a morphism $\Phi : \DBDual{X} \to
\omega_X^{\mydot}$, which is dual to the natural morphism $\sO_X \to \DuBois{X}$.
\begin{remark}
\label{rem.EnoughToProveDualQuasiIsomorphism}
Note that $X$ has DB singularities if and only if $\Phi$ is a quasi-isomorphism since
applying the Grothendieck duality functor again yields a morphism $\O_X \to
\DuBois{X}$ which can be identified with the morphism from Lemma
\ref{lem.BasicPropertiesOfDuBois}(ii) up to quasi-isomorphism.
\end{remark}

\section{The key injectivity}

\noindent
In Theorem \ref{thm.MainInjectivity} below, we prove the
following injectivity.  For every integer $j \in \bZ$,
\begin{equation*}
\Phi^j : \myH^j(\DBDual{X}) \into \myH^j(\omega_X^{\mydot})
\end{equation*}
is injective.  In the case that $x \in X$ is a closed point such that $X \setminus \{
x \}$ is DB, the injectivity of this morphism played a key role in proving that rational,
log canonical and $F$-injective singularities are DB, see \cite{KovacsDuBoisLC1, KollarKovacsLCImpliesDB,SchwedeFInjectiveAreDuBois}.

Because of its potential usefulness it has been asked several times whether this
injectivity holds.  In particular, it was asked in \cite[Question
8.3]{SchwedeFInjectiveAreDuBois} and \cite[Question 5.2]{KovacsSchwedeDuBoisSurvey}.

First we prove a lemma that is interesting on its own.

\begin{lemma}
\label{lem.CompatibilityOfCyclicCoverAndDuBois}
Let $X$ be a variety and $\sL$ a semi-ample line bundle.  Choose $s \in \sL^n$ a
general global section for some $n\gg 0$ and take the $n^{\text{th}}$-root of this
section as in \cite[2.50]{KollarMori}:
$$
\eta : Y = \sheafspec \bigoplus_{i = 0}^{n-1} \sL^{-i} \to X.
$$
Then $\eta_*=\myR\eta_*$,
\[
\eta_* \DuBois{Y} \simeq \DuBois{X}\otimes \eta_*\sO_Y \simeq \bigoplus_{i = 0}^{n-1}
(\DuBois{X} \tensor \sL^{-i}),
\]
and this direct sum decomposition is compatible with the decomposition $\eta_* \O_Y =
\bigoplus_{i = 0}^{n-1} \sL^{-i}$.
\end{lemma}

\begin{proof}
  We fix $\pi_{\mydot} : X_{\mydot} \to X$ a finite cubic (or simplicial)
  hyperresolution of $X$ as in \cite{GNPP}.  On each component $X_i$ of $X_{\mydot}$,
  $\sL$ pulls back to a semi-ample line bundle and further $s$ is still a general
  member of the base-point free linear subsystem of $\pi_i^* \sL^n$.  Thus we obtain
  a cyclic cover $\eta_i : Y_i \to X_i$ for each $i$ as well.  Furthermore, each
  $Y_i$ is smooth since it is ramified over a general element of a base-point free
  linear system.  Obviously, these $Y_i$'s glue to give a diagram of smooth
  $\bC$-schemes $Y_{\mydot}$ with an augmentation morphism $\rho_{\mydot} :
  Y_{\mydot} \to Y$.  From the construction of a cubic hyperresolution, it is easy to
  see that $Y_{\mydot}$ is also a cubic hyperresolution.

  We briefly sketch the idea of this last claim: if $X' \to X$ is a resolution of
  singularities, then the induced $Y' \to Y$ is also a resolution of singularities.
  Furthermore, if $X' \to X$ is an isomorphism outside of $\Sigma \subseteq X$, then
  $Y' \to Y$ is also an isomorphism outside of $\eta^{-1}(\Sigma)$, which is itself
  the induced cyclic cover of $\Sigma$.

Therefore,
\begin{align*}
  \myR\eta_* \DuBois{Y} & 
  \simeq \myR\eta_* \myR {\rho_{\mydot}}_* \O_{Y_{\mydot}}
  \simeq \myR {\pi_{\mydot}}_* {\myR\eta_\mydot}_*\O_{Y_{\mydot}} \\
  & \simeq \myR{\pi_{\mydot}}_* \left(\oplus_{i=0}^{n-1} (\O_{X_{\mydot}} \tensor
    \pi_{\mydot}^* \sL^{-i})\right) \\
  & \simeq \oplus_{i = 0}^{n-1} \left( \left(\myR{\pi_{\mydot}}_*
      \O_{X_{\mydot}}\right) \tensor \sL^{-i} \right) \\
  & \simeq \oplus_{i = 0}^{n-1} (\DuBois{X} \tensor \sL^{-i})\\
  & \simeq \DuBois{X} \tensor (\oplus_{i = 0}^{n-1}  \sL^{-i}) \\
  & \simeq \DuBois{X}\otimes \eta_*\sO_Y.
\end{align*}
and the result follows, the compatibility statement following by construction.

Alternatively, if one wishes to avoid hyperresolutions one may proceed as follows.
By restricting to an open set, we may assume that $X$ embeds as a closed subscheme in
a smooth scheme $U$ such that $\sL$ is the restriction of a globally generated
line-bundle $\sM$ on $U$.  Further set $\pi : U' \to U$ to be a log resolution of $X
\subseteq U$ where we use $\overline{X}$ to denote the reduced divisor
$\pi^{-1}(X)_{\reduced}$.  Then $\myR \pi_* \O_{\overline X} \simeq \DuBois{X}$.
Choosing a general section $s$ of the globally generated line bundle $\sM^{n}$, we
obtain a diagram of cyclic covers:
\[
\xymatrix{
  \overline{Y} \ar@{^{(}->}[r] \ar[d] & W' \ar[d] \\
  Y \ar@{^{(}->}[r] & W }
\]
where $Y, W, W'$ and $\overline{Y}$ are the induced cyclic covers of $X, U, U'$ and
$\overline{X}$ respectively.  It is clear that $W$ and $W'$ are smooth and that
$\overline{Y}$ is the reduced-preimage of $Y$ and has simple normal crossings.  Thus
the result follows again since $\myR \pi_* \sO_{\overline{Y}} \simeq \DuBois{Y}$ by
\cite{SchwedeEasyCharacterization}, also see
\cite{EsnaultHodgeTypeOfSubvarietiesOfPn}.
\end{proof}

Before proving our main injectivity, we need one more result.

\begin{proposition}
\label{prop.SurjectivityWithSemiample}
Let $X$ be a proper variety over $\bC$ and $\sL$ a semi-ample line bundle on $X$.
Then the natural map
\[
H^j(X, 
\sL^{-i}) \to \bH^j(X, \DuBois{X} \tensor \sL^{-i})
\]
is surjective for all $j, i \geq 0$.
\end{proposition}
\begin{proof}
  Choose $n > i$ such that $\sL^{n}$ is base-point-free and choose a general section
  $s \in \Gamma(X, \sL^{n})$. Consider the induced cyclic cover $\eta : Y \to X$ and
  note that $Y$ is also proper.  Now, we have the following factorization
  \[
  H^i(Y^{\an}, \bC) \to H^i(Y, \O_Y) \to \bH^i(Y, \DuBois{Y}).
  \]
  This composition is surjective by Lemma \ref{lem.BasicPropertiesOfDuBois}(iii).
  Thus $H^i(Y, \O_Y) \to \bH^i(Y, \DuBois{Y})$ is also surjective. Then the statement
  follows by Lemma \ref{lem.CompatibilityOfCyclicCoverAndDuBois}.
\end{proof}

Now we are ready to prove the main result of the section.

\begin{theorem}
  \label{thm.MainInjectivity}
  Let $X$ be a variety over $\bC$.  Then the natural map
  \[
  \Phi^j : \myH^j(\DBDual{X}) \into \myH^j(\omega_X^{\mydot})
  \]
  is injective for every $j\in\bZ$.
\end{theorem}
\begin{proof}
  The statement is local and compatible with restriction to an open subset. Therefore
  we may assume that $X$ is projective.  Let $j\in\bZ$ and $\sL$ an ample line bundle
  on $X$.  It follows from Proposition \ref{prop.SurjectivityWithSemiample} that
  $H^{-j}(X, \sL^{-i}) \to \bH^{-j}(X, \DuBois{X} \tensor \sL^{-i})$ is surjective.
  Next, apply $\Hom_{\bC}(\blank, \bC)$ and observe that then
  \[
  \bH^{-j}(X, \DuBois{X} \tensor \sL^{-i})^{\vee} \into H^{-j}(X, \sL^{-i})^{\vee}
  \]
  is injective.  However,
  $$
  H^{-j}(X, \sL^{-i})^{\vee} \simeq \myH^{j}(\myR \Gamma(X, \myR \sHom_{\O_X}(\sL^{-i},
  \omega_X^{\mydot} ))) \simeq \bH^{j}(X, \omega_X^{\mydot} \tensor \sL^{i})
  $$ by
  Grothendieck duality applied to the structure map $\epsilon:X\to \bC$.  Likewise,
  $$
  \bH^{-j}(X, \DuBois{X} \tensor \sL^{-i})^{\vee} \simeq \bH^{j}(X, \DBDual{X} \tensor
  \sL^{i}).
  $$
  Thus we get that
  \begin{equation*}
    \bH^{j}(X, \DBDual{X} \tensor \sL^{i}) \into \bH^{j}(X, \omega_X^{\mydot} \tensor \sL^{i})
  \end{equation*}
  is injective.  Notice that $\bH^{j}(X, \omega_X^{\mydot} \tensor \sL^{i}) \simeq
  H^{0}(X, \myH^{j}(\omega_X^{\mydot}) \tensor \sL^{i})$ for $i \gg 0$ by
  Serre-vanishing and the associated Grothendieck spectral sequence.  Likewise,
  $\bH^{j}(X, \omega_X^{\mydot} \tensor \sL^{i}) \simeq H^{0}(X,
  \myH^{j}(\DBDual{X}) \tensor \sL^{i})$ for $i \gg 0$.  Therefore,
  \begin{equation}
  \label{eq.InjectivityAfterSS}
  H^{0}(X,  \myH^{j}(\DBDual{X}) \tensor \sL^{i}) \into
  H^{0}(X, \myH^{j}(\omega_X^{\mydot}) \tensor \sL^{i})
  \end{equation}
  is injective for $i \gg 0$.  Observe that since $\sL$ is ample, both
  $\myH^{j}(\DBDual X)\tensor \sL^{i}$ and $\myH^{j}(\omega_X^{\mydot})\tensor
  \sL^{i}$ are generated by global sections for $i\gg 0$. Therefore the injectivity
  of equation (\ref{eq.InjectivityAfterSS}) implies, that
  \[
  \Phi^j : \myH^j(\DBDual{X}) \to \myH^j(\omega_X^{\mydot})
  \]
  is also injective for every $j$.  This completes the proof.
\end{proof}

We also have the following local-dual version of Theorem \ref{thm.MainInjectivity}.

\begin{corollary}[(\protect{\cf~\cite[Lemma 2.2]{KovacsDuBoisLC1}})]
  \label{cor.LocalDualSurjectivity}
  Let $X$ be a variety and $P \in X$ is a point (not necessarily closed).  Then the
  natural map
  \[
  H^i_{P}( X, \O_{X, P}) \onto \bH^i_P(X, \DuBois{X}\tensor \O_{X,P})
  \]
  is surjective for all $i \geq 0$.
\end{corollary}
\begin{proof}
  We have the injection $\myH^i (\DBDual{X})_P \to \myH^i (\omega_X^{\mydot})_P$ for
  all $i$.  After shifting (in case $P$ is not a closed point), we have that $\myH^i
  (\DBDual{\O_{X,P}}) \to \myH^i(\omega_{\O_{X, P}}^{\mydot})$ also injects for all
  $i$. Let $E$ be the injective hull of the residue field $\O_{X, P}/\frm_{X,P}$ and apply the
  (faithful and exact) functor $\Hom_{\O_{X, P}}(\blank, E)$.  Local duality in the
  form of \cite[IV, Theorem 6.2] {HartshorneResidues} then yields the corollary.
\end{proof}

With respect for deciding whether $X$ has DB singularities, the complex $\DuBois{X}$ plays the same role as the complex $\myR \pi_*
\O_{\tld X}$ does for detecting rational singularities, here $\pi : \tld X \to X$ is a resolution of singularities.

However, in many applications what makes $\myR \pi_*
\O_{\tld X}$ a useful object is Grauert-Riemenschneider vanishing \cite{GRVanishing} applied to it's Grothendieck dual.  In particular, the
Grothendieck dual $\myR \pi_* \omega_{\tld X}^{\mydot} \qis \myR
\sHom_{\O_{X}}^{\mydot}( \myR \pi_* \O_{\tld X}, \omega_X^{\mydot})$ is a complex
with cohomology in only one spot, $$\myR \pi_* \omega_{\tld X}^{\mydot} \qis \pi_*
\omega_{\tld X}[\dim X].$$

For $X$ Cohen-Macaulay, Theorem \ref{thm.MainInjectivity} yields an analogous
vanishing for DB singularities.

\begin{corollary}
  \label{cor.GRVanishingForDuBois}
  Let $X$ be a Cohen-Macaulay variety of dimension $d$.  Then
  \[
  \DBDual{X} \qis \myH^{-d} (\DBDual{X}) [d].
  \]
  If additionally $X$ is normal and $\pi : \tld X \to X$ is a log resolution of
  singularities with reduced exceptional divisor $E$, then
  \[
  \DBDual{X} \qis \pi_* \omega_{\tld X}(E)[d]
  \]
\end{corollary}
\begin{proof}
  Since $X$ is Cohen-Macaulay and connected, it is equidimensional.  The first
  statement is immediate since a submodule of the zero-module is zero and because
  $\myH^i(\omega_X^{\mydot}) = 0$ for $i \neq -d$.  For the second statement, use the
  fact that $\myH^{-d} (\DBDual{X}) \simeq \pi_* \omega_{\tld X}(E)$ by
  \cite[Theorem 3.8]{KovacsSchwedeSmithLCImpliesDuBois}.
\end{proof}

\begin{remark}
  Notice that if $X$ is DB, then $\DBDual{X}\simeq \omega_X^\mydot$ and hence the
  statement is equivalent to $X$ being Cohen-Macaulay.
\end{remark}

A slight reinterpretation of the previous result also gives us the following corollary.
\begin{corollary}
  Let $Y$ be a smooth $n$-dimensional variety and $X \subseteq Y$ a Cohen-Macaulay
  subvariety of pure dimension $d$.  Let $\pi : \tld Y \to Y$ be a log resolution of
  $X \subseteq Y$.  Set $E \subseteq Y$ to be the reduced pre-image of $X$ in $Y$
  (which is a divisor since $\pi$ is a log resolution).  Then
  \[
  R^i \pi_* \omega_{\tld Y}(E) = 0
  \]
  for all $i \neq 0, n - d - 1$.
\end{corollary}
\begin{proof}
  Consider the long exact sequence
  \[
  R^i \pi_* \omega_{\tld Y} \to R^i \pi_* \omega_{\tld Y}(E) \to R^i \pi_* \omega_{E}
  \to R^{i+1} \pi_* \omega_{\tld Y}
  \]
  and notice first that $R^i \pi_* \omega_{\tld Y} = 0$ for all $i \neq 0$ by
  \cite{GRVanishing}.  Since $\omega_E[n-1] \simeq \omega_E^{\mydot}$ we have $R^{j +
    n - 1} \pi_* \omega_E \qis \myH^{j}(\myR \pi_* \omega_E^{\mydot})$.  However,
  $\myR \pi_* \omega_E^{\mydot} \qis \DBDual{X}$ by
  \cite{SchwedeEasyCharacterization}.  Therefore, since $\myH^j(\DBDual{X}) = 0$ for
  $j \neq -d$ by Corollary \ref{cor.GRVanishingForDuBois}, we see that $R^{j + n - 1}
  \pi_* \omega_E = 0$ for $j \neq -d$.  Thus $R^i \pi_* \omega_E = 0$ for $i \neq n -
  d - 1$ and the result follows.
\end{proof}

\begin{remark}
  The previous two corollaries do not hold if $X$ is not Cohen-Macaulay.  In fact
  they automatically fail for any non-Cohen-Macaulay variety with Du Bois
  singularities.  For example, they fail for the affine cone over an Abelian variety
  of dimension $> 1$.
\end{remark}

Theorem \ref{thm.MainInjectivity} also provides slightly simpler proofs of existing results.

\begin{corollary}[(\cite{KovacsDuBoisLC1}, \protect{\cf~\cite[Section
    12]{KollarShafarevich}})] If the morphism $\O_X \to \DuBois{X}$ has a
  left-inverse in $D^{b}_{\coherent}(X)$, then $X$ has DB singularities.
\end{corollary}

\begin{proof}
  The hypothesis implies that $\Phi^i : \myH^i(\DBDual{X}) \to \myH^i
  (\omega_X^{\mydot})$ is surjective for every $i$.  Thus $\Phi^i$ is an isomorphism
  by Theorem \ref{thm.MainInjectivity} and hence $\Phi : \DBDual{X} \to
  \omega_X^{\mydot}$ is a quasi-isomorphism and so $X$ has DB singularities by Remark
  \ref{rem.EnoughToProveDualQuasiIsomorphism}.
\end{proof}

\ifthenelse {\boolean{bool.doWithPairs}} { }
{
\begin{remark}
In \cite{KovacsDBPairsAndVanishing}, the first author of this paper introduced the notion of Du~Bois pairs.  It is easy to see that in the language of that paper, that an analog of Theorem \ref{thm.MainInjectivity} holds.
In particular, if $(X, \Sigma)$ is a pair, then the natural map
\[
  \Phi^j : \myH^j( \myR \sHom_{\O_X}(\DuBois{X,\Sigma}, \omega_X^{\mydot})) \into \myH^j(\myR \sHom_{\O_X}(\sI_{\Sigma}, \omega_X^{\mydot}))
\]
is injective for ever $j$.
\end{remark}
}

\section{Deformation of DB singularities}

We now prove the main result of the paper.  In fact, simply using Corollary
\ref{cor.LocalDualSurjectivity} fills in the gap in the first author's proof of this
statement in \cite[Theorem 3.2]{KovacsDuBoisLC2}.  For completeness, we provide a
proof below.  This proof (as well as the proof of \cite[Theorem
3.2]{KovacsDuBoisLC2}) was inspired by Elkik's proof of the fact that rational
singularities deform \cite{ElkikDeformationsOfRational}.

\begin{theorem}
  \label{thm.MainTheorem}
  Let $X$ be a scheme of finite type over $\bC$ and $H$ a reduced effective Cartier
  divisor (if $X$ is not normal, by a Cartier divisor we mean a subscheme locally
  defined by a single non-zero-divisor at each stalk).  If $H$ has DB singularities,
  then $X$ has DB singularities near $H$.
\end{theorem}
\begin{proof}
  Choose hyperresolutions $\pi_{\mydot} : X_{\mydot} \to X$ and $\mu_{\mydot} :
  H_{\mydot} \to H$ with a map $H_{\mydot} \to X_{\mydot}$ factoring through the
  diagram of schemes $Z_{\mydot} := X_{\mydot} \times_X H$ as pictured below, \cf
  \cite{GNPP}.
  \[
  \xymatrix{ H_{\mydot} \ar[rd]_{\mu_{\mydot}} \ar[r] & Z_{\mydot}
    \ar[d]^-{\varepsilon_{\mydot}} \ar[r] & X_{\mydot} \ar[d]^-{\pi_{\mydot}} \\
    & H \ar@{^{(}->}[r] & X }
  \]
  Note that the components of $Z_{\mydot}$ need not be smooth or even reduced.

  Choose a closed point ${\bq}$ of $X$ contained within $H$.  It is sufficient to
  prove that $X$ is DB at ${\bq}$.  Let $R$ denote the stalk $\O_{X, {\bq}}$ and
  choose $f \in R$ to denote a defining equation of $H$ in $R$.  We also define
  $\DuBois{R} := \DuBois{X} \tensor R$ and $\DBDual{R} := \myR
  \Hom_R^{\mydot}(\DuBois{R}, \omega_R^{\mydot})$.  Consider the following diagram
  whose rows are exact triangles in $D^{b}_{\coherent}(X)$:

  \begin{equation*}
    \xymatrix@R=12pt{
      R \ar[d] \ar[r]^{\times f} & R \ar[d] \ar[r] & R/(f) \ar[d]_{\rho} \ar[r]^-{+1}
      & \\
      \DuBois{R} \ar[r]_{\times f} & \DuBois{R} \ar[r] & \left(\myR
      \varepsilon_{\mydot *} \sO_{Z_{\mydot}}\right) \tensor  R
      \ar[d]_{\tau}      \ar[r]^-{+1} & \\
      & & \DuBois{H}\tensor R
      &\\
    }
  \end{equation*}
  where $\tau \circ \rho$ is a quasi-isomorphism by hypothesis.
  Next we apply the functor $\myR \Hom_R^{\mydot}(\blank, \omega_R^{\mydot})$.
  Using the notation ${\wt\omega}_{Z_\mydot}^{\mydot}=\myR \Hom_R^{\mydot}(
  \left(\myR \varepsilon_{\mydot *}\sO_{Z_{\mydot}}\right) \tensor R
  , \omega_R^{\mydot})$ and taking cohomology we obtain the diagram of long exact
  sequences:
  \[
  \xymatrix@R=15pt{%
    \dots \ar@{<-}[r] & \myH^i(\omega_R^{\mydot}) \ar@{<-_{)}}[d]^{\Phi^i} &
    \ar[l]_{\times f} \myH^i(\omega_R^{\mydot}) \ar@{<-_{)}}^{\Phi^i}[d]
    \ar@{<-}[r]^{\delta_i} & \myH^i(\omega_{R/f}^{\mydot}) \ar@{<<-}[d]^-{\gamma_i}
    \ar@{<-}[r]^-{\alpha_i} & \myH^{i-1}(\omega_R^{\mydot})
    \ar@{<-_{)}}[d]^{\Phi^{i-1}} \ar@{<-}[r]^{\times f} &
    \myH^{i-1}(\omega_R^{\mydot}) \ar@{<-_{)}}[d]^{\Phi^{i-1}}  & \ar[l]\cdots  \\
    \dots \ar@{<-}[r] & \myH^i(\DBDual{R}) \ar@{<-}[r]_{\times f} & \myH^i
    (\DBDual{R}) \ar@{<-}[r] & \myH^i({\wt\omega}_{Z_\mydot}^{\mydot})
    \ar@{<-}[r]_-{\beta_i} & \myH^{i-1}(\DBDual{R} \ar@{<-}[r]_{\times f}) &
    \myH^{i-1}(\DBDual{R}) }
  \]
  where the vertical $\Phi$ maps are injective because of Theorem
  \ref{thm.MainInjectivity} and the morphism $\gamma_i$ is surjective because $\tau
  \circ \rho$ is an isomorphism.

  Fix $z \in \myH^{i-1} (\omega_R^{\mydot})$.  Pick $w \in
  \myH^i({\wt\omega}_{Z_\mydot}^{\mydot})$ such that $\alpha_i(z) = \gamma_i(w)$.
  Since $\delta_i(\alpha_i(z)) = 0$ and $\Phi^i$ is injective, it follows that there
  exists a $u \in \myH^{i-1} (\DBDual{R})$ such that $\beta_i(u) = w$.  Therefore,
  $\alpha_i(\Phi^{i-1}(u)) = \alpha_i(z)$ and so
  \begin{equation}
    \label{eq.diffInFTimes}
    z - \Phi^{i-1}(u) \in f \cdot \myH^{i-1} (\omega_R^{\mydot}).
  \end{equation}

  Now, fix $C_{i-1}$ to be the cokernel of $\Phi^{i-1}$ and set $\overline{z} \in
  C_{i-1}$ to be the image of $z$.  Equation (\ref{eq.diffInFTimes}) then guarantees
  that $\overline{z} \in f \cdot C_{i-1}$.  But $z$ was arbitrary and so the
  multiplication map $\xymatrix{C_{i-1} \ar[r]^{\times f} & C_{i-1}}$ is surjective.
  But this contradicts Nakayama's lemma unless $C_{i-1} = 0$.  Therefore $C_{i-1} =
  0$ and $\Phi^{i-1}$ is also surjective.  This holds for all $i$ and so the natural
  morphism $\DBDual{X} \to \omega_X^{\mydot}$ is a quasi-isomorphism.  Thus $X$ has
  DB singularities by Remark \ref{rem.EnoughToProveDualQuasiIsomorphism}.
\end{proof}

\begin{corollary}
  \label{cor.FamilyDeforms}
  Let $f : X \to S$ be a proper flat family of varieties over a smooth curve $S$ and $s \in S$ a
  closed point.  If the fiber $X_s$ has DB singularities, then so do the other fibers
  near $s$.
\end{corollary}
\begin{proof}
  By Theorem \ref{thm.MainTheorem}, $X$ has DB singularities near $X_s$.  Let $\Sigma$ denote the non-Du Bois locus of $X$.  Since $f$ is proper, $f(\Sigma)$ is a closed subset of $S$ not containing $s \in S$.  Thus by
  restricting $S$ to an open set, we may assume that $X$ has DB singularities.  By
  Lemma \ref{lem.GeneralSectionOfDuBois}, all fibers over nearby points of $s \in S$ have
  DB singularities.
\end{proof}

\ifthenelse {\boolean{bool.doWithPairs}} {
  \section{DB pairs}

  In \cite{KovacsDBPairsAndVanishing}, the first author defined a notion of Du~Bois
  (or simply DB) pairs. Indeed, given a (possibly non-reduced) subscheme $Z \subseteq
  X$ one has an induced map in $D^b_{\coherent}(X)$,
  \[
  \DuBois{X} \to \DuBois{Z},
  \]
  noting that by definition $\DuBois{Z} = \DuBois{Z_{\reduced}}$.  Then $\DuBois{X,
    Z}$ to be the object in the derived category making the following an exact
  triangle:
  \[
  \DuBois{X, Z} \to \DuBois{X} \to \DuBois{Z} \xrightarrow{+1}.
  \]
  If $\sI_Z$ is the ideal sheaf of $Z$, then it is easy to see that there is a
  natural map $\sI_Z \to \DuBois{X,Z}$, \cite[Section
  3.D]{KovacsDBPairsAndVanishing}.

  \begin{definition}\cite[Definition 3.13]{KovacsDBPairsAndVanishing}
    The \emph{Du Bois defect of $(X, Z)$}, denoted $\XDuBois{X,Z}$, is the mapping
    cone of the morphism $\sI_Z \to \DuBois{X,Z}$, so that there is an exact triangle
    \[
    \sI_Z \to \DuBois{X,Z} \to \XDuBois{X,Z} \xrightarrow{+1}
    \]
    We say that $(X, Z)$ has \emph{Du Bois singularities} if $\XDuBois{X,Z}$ is
    quasi-isomorphic to zero.  In other words, if $\sI_Z \to \DuBois{X,Z}$ is a
    quasi-isomorphism.
  \end{definition}

  We now mimic our approach before:

  \begin{lemma} [\cf Lemma \ref{lem.CompatibilityOfCyclicCoverAndDuBois}]
    \label{lem.CompatibilityOfCyclicCoverAndDuBoisPairs}
    Let $X$ be a variety and $\sL$ a semi-ample line bundle.  Choose $s \in \sL^n$ a
    general global section for some $n\gg 0$ and take the $n^{\text{th}}$-root of
    this section as in \cite[2.50]{KollarMori}:
    $$
    \eta : Y = \sheafspec \bigoplus_{i = 0}^{n-1} \sL^{-i} \to X.
    $$
    Set $W = \eta^{-1}(Z)$ (with the induced scheme structure).
    Note that we have $\eta|_W : W = \sheafspec \bigoplus_{i = 0}^{n-1} \sL^{-i}|_Z \to Z$.
    Then as before $\eta_*=\myR\eta_*$,
    \[
    \eta_* \DuBois{Y,W} \simeq \DuBois{X,Z}\otimes \eta_*\sO_Y \simeq \bigoplus_{i =
      0}^{n-1} (\DuBois{X,Z} \tensor \sL^{-i}),
    \]
    and this direct sum decomposition is compatible with the decomposition $\eta_*
    \O_Y = \bigoplus_{i = 0}^{n-1} \sL^{-i}$.
  \end{lemma}
  \begin{proof}
    This can be proven just as in Lemma \ref{lem.CompatibilityOfCyclicCoverAndDuBois}
    or alternately follows formally from Lemma
    \ref{lem.CompatibilityOfCyclicCoverAndDuBois} via the functoriality of the
    construction.
  \end{proof}

  Just as in Proposition \ref{prop.SurjectivityWithSemiample}, we also obtain that
  \[
  H^j(X, \sI_Z \tensor \sL^{-i}) \to \bH^j(X, \DuBois{X, Z} \tensor \sL^{-i})
  \]
  is simply using \cite[Theorem 4.1]{KovacsDBPairsAndVanishing} in place of Lemma
  \ref{lem.BasicPropertiesOfDuBois}(iii).

  If we set $\DBDual{X,Z} = \myR \sHom_{\O_X}^{\mydot}(X, \omega_X^{\mydot})$, then
  we easily obtain.

  \begin{theorem}
    \label{thm.MainInjectivityForPairs}
    Let $X$ be a variety over $\bC$.  Then the natural map
    \[
    \Phi^j : \myH^j(\DBDual{X}) \into \myH^j(\myR \sHom_{\O_X}(\sI_Z,
    \omega_X^{\mydot}))
    \]
    is injective for every $j\in\bZ$.
  \end{theorem}
  \begin{proof}
    The proof is the same as in Theorem \ref{thm.MainInjectivity}
  \end{proof}

\newpage

\input tranversality.tex

\newpage

} { }

\section{Application to restriction theorems for maximal non-LC ideals}

In this section we assume the reader is familiar with log canonical singularities;
see \cite{KollarMori} for an introduction.  Let $X$ be a normal variety, $\Delta$ an
effective $\bQ$-divisor on $X$ such that $K_X + \Delta$ is $\bQ$-Cartier and $\pi :
\tld X \to X$ is a log resolution for $(X, \Delta)$.  Write $K_{\tld X} - \pi^* (K_X
+ \Delta) = \sum a_i E_i$ and set $E^{=-1} = \sum_{a_i = -1} E_i$.  The following
ideal
\[
\mJ_{\textnormal{NLC}}(X, \Delta) := \pi_* \O_{\tld X}(\lceil K_{\tld X} - \pi^* (K_X
+ \Delta) + E^{=-1} \rceil)
\]
is defined to be the \emph{non-log canonical ideal of $X$}.

This ideal was first defined by F.~Ambro in \cite[Definition 4.1]{AmbroQuasiLog}
where it was denoted by $\mathcal{I}_{X_{-\infty}}$.  The study of this ideal as an
object similar to the multiplier ideal, was recently initiated by O.~Fujino in
\cite{FujinoNonLCSheaves}.  One of the main facts about this ideal is that the zero
set of $\mJ_{\textnormal{NLC}}$ is exactly the locus where $(X, \Delta)$ does not
have log canonical singularities. Fujino proved the following restriction theorem for
$\mJ_{\textnormal{NLC}}(X, \Delta)$ (in fact, he proved a more general result):

\begin{theorem*} \cite[Theorem 1.2]{FujinoNonLCSheaves} If $H$ is a normal Cartier
  divisor on a $\bQ$-Gorenstein variety $X$, then $\mJ_{\textnormal{NLC}}(X, H)
  \tensor \O_H \simeq \mJ_{\textnormal{NLC}}(H, 0)$.
\end{theorem*}

However, there are other natural ideals that define the non-lc locus.  With notation
as above, set $E = \sum E_i$ and set $E^{\bZ} = \sum_{a_i \in \bZ} E_i$.  Then
consider the ideal
\[
\mJ'(X, \Delta) := \pi_* \O_{\tld X}(\lceil K_{\tld X} - \pi^* (K_X + \Delta) +
E^{\bZ} \rceil) = \pi_* \O_{\tld X}(\lceil K_{\tld X} - \pi^* (K_X + \Delta) +
\varepsilon E \rceil)
\]
where we choose $1 \gg \varepsilon > 0$.  This is the largest ideal which canonically
defines the non-log canonical locus of $(X, \Delta)$ and as such is called the
\emph{maximal non-lc ideal}.  In \cite{FujinoSchwedeTakagiSupplements}, the authors
explored this ideal (and other non-lc-ideals).  In particular, they obtained
restriction theorems in special cases \cite[Theorem 12.7, Theorem
13.13]{FujinoSchwedeTakagiSupplements}.  As an application of Theorem
\ref{thm.MainInjectivity}, we obtain the following restriction theorem for $\mJ'(X,
H)$ in the case that $X$ is Gorenstein.

\begin{theorem}
  \label{thm.RestrictionTheoremForMNLC}
  If $X$ is a normal $d$-dimensional Gorenstein variety and $H$ is a normal Cartier
  divisor on $X$, then $\mJ'(X, H)|_H \simeq \mJ'(H, 0)$.
\end{theorem}
The proof strategy is the same as in \cite[Section
13]{FujinoSchwedeTakagiSupplements}
\begin{proof}
  By working sufficiently locally, we may assume that $K_X \sim 0$ and $H = V(f) \sim
  0$ for some $f \in \Gamma(X, \O_X)$.  Shrinking $X$ again if necessary, we embed $X
  \subseteq Y$ as a closed subscheme in a smooth scheme $Y$.  Let $\pi : \tld Y \to
  Y$ be a log resolution of $H \subseteq Y$ which is simultaneously an embedded
  resolution of $X \subseteq Y$.  Let $\overline{X} = \pi^{-1}(X)_{\reduced}$, $\tld
  X$ the strict transform of $X$, and $\overline{H} = \pi^{-1}(H)_{\reduced}$.  We
  may assume that $\pi$ is an isomorphism outside of $\Sing X \cup H$ and write
  $\overline{X} = \tld{X} \cup E \cup \overline{H}$ where $E = \pi^{-1}(\Sing
  X)_{\red}$.  Finally, we may also assume that $E \cup \overline H$ is a reduced
  simple normal crossings divisor which intersects $\tld X$ with normal crossings so
  that $(E \cup \overline{H}) \cap \tld{X}$ is a reduced simple normal crossings
  divisor on $\tld X$.  We have the following short exact sequence:
  \[
  0 \to \O_{\tld X}(-E \cup \overline{H}) \to \O_{\overline{X}} \to \O_{E \cup
    \overline{H}} \to 0.
  \]
  By pushing forward and using \cite{SchwedeEasyCharacterization}, we obtain the
  exact triangle,
  \[
  \xymatrix{ \myR \pi_* \O_{\tld X}(-E \cup \overline{H}) \ar[r] & \DuBois{X} \ar[r]
    & \DuBois{H \cup \Sing X} \ar[r]^-{+1} & }.
  \]
  Applying $\myR \sHom_{\O_X}^{\mydot}(\blank, \omega_X^{\mydot})$ gives
  \[
  \xymatrix{ \DBDual{H \cup \Sing X} \ar[r] & \DBDual{X} \ar[r] & \myR \pi_* \O_{\tld
      X}(K_{\tld X} + E \cup \overline{H})[d] \ar[r]^-{+1} & },
  \]
  and by taking cohomology, we arrive at the exact sequence
  {\small
  \begin{equation}
  \label{eq.ExactSequenceForDBDualsAdjunction}
  0 \to \myH^{-d} (\DBDual{X}) \to \pi_* \O_{\tld X}(K_{\tld X} + E \cup
  \overline{H}) \to \myH^{-d + 1} (\DBDual{H \cup \Sing X}) \to \myH^{-d+1}
  (\DBDual{X}) = 0.
  \end{equation}}\noindent
The vanishing on the right follows by Corollary \ref{cor.GRVanishingForDuBois}
  since $X$ is Gorenstein and thus Cohen-Macaulay.

  By \cite[Lemma 13.11]{FujinoSchwedeTakagiSupplements}, $\myH^{-d + 1} (\DBDual{H
    \cup \Sing X}) \simeq \myH^{-d + 1} (\DBDual{H})$.  Furthermore, by \cite[Theorem 3.8]{KovacsSchwedeSmithLCImpliesDuBois} we know $\mJ'(X, 0) \cong h^{-d}(\DBDual{X}) \tensor \O_X(-K_X)$ and $\mJ'(H, 0) \cong h^{-d+1}(\DBDual{H}) \tensor \O_X(-K_X - H)$.  Hence twisting \eqref{eq.ExactSequenceForDBDualsAdjunction} by $\O_X(-K_X-H)$ we obtain the following short exact sequence: \cf \cite[Lemma
  13.8]{FujinoSchwedeTakagiSupplements}  \cite[Lemma
  4.14]{KovacsSchwedeSmithLCImpliesDuBois},
  \[
  0 \to \mJ'(X, 0) \tensor \O_X(-H) \to \pi_* \O_{\tld X}(K_{\tld X} - \pi^*(K_X + H) + E \cup \overline{H})
  \to \mJ'(H, 0) \to 0.
  \]
This completes the proof.
\end{proof}

\bibliographystyle{skalpha}
\bibliography{CommonBib}

\end{document}